\documentstyle[12pt]{article}
\newtheorem{theorem}{THEOREM}[section]
\newtheorem{definition}[theorem]{Definition}
\newtheorem{corollary}[theorem]{Corollary}
\newtheorem{lemma}[theorem]{Lemma}
\newtheorem{example}{EXAMPLE}
\newtheorem{remark}{REMARK}
\newtheorem{proposition}[theorem]{Proposition}
\newtheorem{problem}{PROBLEM}
\def\qbar{\overline{q}}
\def\lbar{\overline{\ell}}

\def\jbar{\overline{j}}
\def\kbar{\overline{k}}

\def\L{{\cal L}}
\def\hbar{\overline{h}}

\def\CC{{\rm\kern.24em\vrule
width.02em height1.4ex depth-.05ex\kern-.26em C}}
\def\QQ{{\rm\kern.24em\vrule width.02em
height1.4ex depth-.05ex\kern-.26em Q}}
\def\RR{{\rm I\kern-.2em R}}
\def\HH{{\rm I\kern-.2em H}}
\def\ZZ{{\rm\kern.26em\vrule width.02em
height0.5ex depth0ex\kern.04em\vrule width.02em height1.47ex
depth-1ex\kern-.34em Z}}
\def\Ibb#1{{\rm I\kern-.23em#1}}
\def\Ib#1{{\rm I\kern-.25em#1}}
\def\k#1{\kern#1em}
\def\vb#1{\vrule width.02em height1.4ex depth-.05ex}

\def\11{{\rm\k{.45}\vb0\k{-.142}1}}

\def\Re{\mbox{\rm Re}\,}

\def\epf{\hskip.2in\vrule width.4pt height6.65pt
depth.15pt\vrule width2.5pt height6.65pt
depth-6.25pt\hskip-2.5pt\vrule width2.5pt height.25pt
depth.15pt\vrule width.4pt height6.65pt depth.15pt\ }
\def\proof{\noindent {\bf Proof. }}

\def \tr{\hbox{tr}}

\def\zbar{\overline{z}}
\def\ibar{\overline{i}}
\def\jbar{\overline{j}}
\def\dbar{\overline{\partial}}
\def \hbar{\overline{h}}

\def\d{\partial}
\def\vbar{\overline{v}}
\def\1bar{\overline{1}}

\def \Dbar{\overline{D}}

\def\nbar{\overline{n}}

\def\Dbar{\overline{D}}

\def\tbar{\overline{t}}

\def\TN{\tilde{\nabla}}
\def\Rbar{\overline{R}}
\def\TD{\tilde{\Delta}}

\begin{document}

\title{ On plurisubharmonicity of the solution of the Fefferman equation and its applications to estimate the  bottom of the spectrum of Laplace-Beltrami operators}

\author{Song-Ying Li}

\date{Revised January 28, 2015}
\maketitle

{ \noindent {\bf Abstract:}
In this paper, we introduce a concept of  super-pseudoconvex domain. We prove that the solution 
of the Feffereman equation on a smoothly bounded strictly pseudoconvex domain $D$ in $\CC^n$ is  plurisubharmonic
if and only if $D$ is super-pseudoconvex. As an application, 
we give a lower bound  estimate the bottom of the spectrum of Laplace-Beltrami operators when $D$ is super-pseudoconvex
by using the result of Li and Wang \cite{LiWang}.
}

\section{Introduction}

Let  $D$ be a smoothly bounded pseudoconvex domain $D$  in $\CC^n$. Let $u\in C^2(D)$ be a real-valued function and let $H(u)$ denote
the $n\times n$ complex Hessian matrix of $u$. We say that $u$ is strictly plurisubharmonic in $D$ if $H(u)$ is positive definite
on $D$. When $u$ is strictly plurisubharmonic in $D$, $u$ induces a K\"ahler metric
$$
g=g[u] =\sum_{i,j=1}^n {\d^2 u\over \d z_i \d \zbar_j} dz^i \otimes d \zbar^j.\leqno(1.1)
$$
We say that the metric $g$ is also Einstein if its Ricci curvature
$$
R_{k\lbar}=-{\d^2 \log \det [g_{i\jbar}] \over \d z_k \d \zbar_{\ell}}=c g_{k\lbar}\leqno(1.2)
$$
for some constant $c$. 

When $c<0$, after a normalization, we may assume $c=-(n+1)$.
It was proved by Cheng and Yau \cite{CY}
that the following Monge-Amp\`ere equation:
$$
\cases{\det H(u)=e^{(n+1) u},   \quad z\in D\cr \quad \qquad u=+\infty,  \qquad z\in \d D\cr}\leqno(1.3)
$$
has a unique strictly plurisubharmonic solution $u\in C^\infty(D)$.  Moreover, the K\"ahler metric
$$
g[u]=\sum_{i,j=1}^n {\d^2 u\over \d z_i \d \zbar_j} dz_i\otimes d\zbar_j\leqno(1.4)
$$
induced by $u$ is a complete K\"ahler-Einstein metric on $D$.

When $D$ is also strictly pseudoconvex,  the existence and uniqueness problem was studied by C. Fefferman \cite{F1} earlier. He considered the
following Fefferman equation
$$
\cases{\det J(\rho)=1,   \quad z\in D\cr \quad \qquad \rho=0,   \quad z\in \d D,\cr}\leqno(1.5)
$$
where
$$
(1.6)\quad J(\rho)=-\det \left[\matrix{ \rho &\dbar \rho \cr (\dbar \rho)^* & H(\rho)\cr}\right],\quad \dbar\rho=({\d \rho\over \d \zbar_1},\cdots, {\d \rho\over \d \zbar_n})
\hbox{ and } (\dbar \rho)^*=({\d \rho\over \d z_1},\cdots, {\d \rho\over \d z_n})^t.
$$
C. Fefferman searched for a solution $\rho<0$ on $D$ such that $u=-\log (-\rho)$ is strictly plurisubharmonic in $D$.  He proved the uniqueness and gave
a formal or approximation solution for (1.5). 

If the relation between $\rho$ and $u$ is given by
$$
\rho(z)=-e^{-u(z)},\qquad z\in D,\leqno(1.7)
$$
then (1.3) is the same as (1.5). Moreover, one can prove (see \cite{Li1} and references therein) that
$$
\det H(u)= J(\rho) e^{(n+1) u}.\leqno(1.8)
$$

When $D$ is smoothly bounded strictly pseudoconvex, it was proved by Cheng and Yau \cite{CY} that $\rho\in C^{n+3/2}(\Dbar)$. In fact,
$\rho\in C^{n+2-\epsilon}(\Dbar)$  for any small $\epsilon>0$. This follows from an asymptotic expansion formula for $\rho$ obtained by Lee and Melrose \cite{LM}:
$$
\rho(z)=r(z)\Big(a_0(z)+\sum_{j=1}^\infty a_j (r^{n+1} \log (-r))^j\Big), \leqno(1.9)
$$
where $r\in C^\infty(\Dbar)$ is any defining function for $D$ and $a_j\in C^\infty(\Dbar)$ and $a_0(z)>0$ on $\d D$. 
\medskip

When $D$ is a bounded strictly pseudoconvex domain in $\CC^n$ with smooth defining function $r$, one can view
$(\d D, \theta)$ as a pseudo-Hermitian CR manifold with the contact/pseudo Hermitian form
$$
\theta={1\over 2 i} (\d r-\dbar r).\leqno(1.10)
$$
An interesting and useful question is:  How to find a defining function $r$ such that $(\d D, \theta)$ has positive
the Webster-Tanaka pseudo Ricci curvature or pseudo scalar curvature? Under the assumption $u=-\log(-r)$ is strictly plurisubharmonic near and on  $\d D$,
the following formula for the pseudo-Ricci curvature was discovered by Li and Luk \cite{LLuk1}:
$$
{\it Ric}_z(w, \vbar)=-\sum_{k,\ell=1}^n {\d ^2\log J(r)(z)\over \d z_k\d \zbar_{\ell}}w_k \vbar_{\ell}+
 n{\det H(r)\over J(r)} \sum_{j,k=1}^n {\d ^2 r(z) \over \d z_k\d \zbar_{\ell}}w_k \vbar_{\ell}\leqno(1.11)
$$
for $w, v\in H_z=\{ v=(v_1,\cdots, v_n)\in \CC^n: \sum_{j=1}^n {\d r(z)\over \d z_j} v_j=0\}$.

When $g[u]$ is asymptotic Einstein (i.e. $J(r)=1+O(r^2))$, one has that
$$
{\it Ric}_z(w, \vbar)=
 n {\det H(r)\over J(r)} \sum_{j,k=1}^n {\d ^2 r(z) \over \d z_k\d \zbar_{\ell}}w_k \vbar_{\ell}\leqno(1.12)
$$
for $w, v\in H_z=\{ v=(v_1,\cdots, v_n)\in \CC^n : \sum_{j=1}^n {\d r(z)\over \d z_j} v_j=0\}$. In this case, the Webster-Tanaka 
pseudo-Hermitian metric is a pseudo Einstein metric. Moreover, it is positive on $\d D$ if and only if $\det H(r)>0$ on
$\d D$. 

Many research works \cite{LT, Li1, Li2, LiWang}  indicate that  the following problem is very interesting and
very important. 

\begin{problem} If  $D$ is a smoothly bounded strictly pseudoconvex domain in $\CC^n$.
Let $\rho$ be the solution of the Fefferman equation (1.5) such that $u=-\log (-\rho)$ is strictly plurisubharmonic in
$D$.  Then $\rho$   is strictly plurisubharmonic in $\Dbar$.
\end{problem}

It is well known that $\rho(z)=|z|^2-1$ is strictly plurisubharmnic when $D=B_n$, the unit ball in $\CC^n$. It was proved by the
Li \cite{Li1} that $\rho$ is strictly plurisubharmonic when $D$ is the bounded domain in $\CC^n$ whose boundary is a real ellipsoid. 
In particular, when $n=2$ case, this result was also proved by Chanillo, Chiu and Yang \cite {CCY} later. \medskip

One of the main purposes of this paper is to give a characterization for domains $D$  in $\CC^n$ where the answer of Problem 1 is affirmatively true.   
We first  introduce the following definition.

\begin{definition} Let $D$ be a smoothly bounded pseudoconvex domain in $\CC^n$. We say that $D$ is strictly super-pseudoconvex (super-superconvex)
 if there is a strictly plurisubharmonic
defining function $r\in C^4(\Dbar)$ such that ${\L}_2[r]>0 $ (${\L}_2[r]\ge 0$) on $\d D$, respectively. Here
$$
{\L}_2[r]=:1+{|\d r|^2_r\over n(n+1)} {\TD} \log J(r)-{2 \Re R \log J(r)\over n+1} -|\d r|^2_r |\TN \log J(r)|^2,\leqno(1.13)
$$
and
$$
\TD =a^{i\jbar}[r]{\d^2\over \d z_i \d \zbar_j},\quad R=\sum_{j=1}^nr^j {\d \over \d z_j}, \quad|\TN f|^2=
a^{i\jbar}[r]{\d f \over \d z_i } {\d f\over \d \zbar_j}\leqno(1.14)
$$
and
$$
r^i =\sum_{j=1}^n r^{i\jbar} r_{\jbar},\quad \Big[r^{i\jbar}\Big]^t=H(r)^{-1},\quad a^{i\jbar}[r]=: r^{i\jbar}-{r^i r^{\jbar}\over -r+|\d r|_r^2},\quad 1\le i,j \le n.\leqno(1.15)
$$
\end{definition}

Another motivation of this paper is to apply the result (the solution of  Problem 1) to estimate the lower bound of the 
bottom of the spectrum of Laplace-Beltrami operator $\Delta_{g[u]}$.

\begin{definition} Let $D$ be a smoothly bounded strictly pseudoconvex domain in $\CC^n$.  Let $r\in C^\infty(\Dbar)$ be
a defining function for $D$ such that $u=-\log (-r)$ is strictly plurisubharmonic. We say that the K\"ahler metric $g[u]$ induced
by $u$ is {\bf super asymptotic Einstein} if 

\item{(i)} the Ricci curvature $R_{i\jbar}\ge -(n+1) g_{i\jbar}$ on $D$; %Equivalently, $-\log J(r)$ is plurisubharmonic;

\noindent and

\item{(ii)} $J(r)=1+O(r^2)$.
\end{definition}
 
 Let $(M^n, g)$ be a K\"ahler manifold with the K\"ahler metric $g$. Let $\Delta_g$ be the Laplcae-Beltrami operator associated to
 $g$. Let $\lambda_1$ denote the bottom of the spectrum of $\Delta_g$. Then estimates of the upper bound and lower
 bound for $\lambda_1$ have studied by many authors, including  S-Y. Cheng\cite{Cheng}, J. Lee \cite{Lee}, P. Li and J-P. Wang \cite{LW1, LW2}, O. Munteanu \cite{M}, S-Y. Li and M-A. Tran \cite{LT} and S-Y. Li and X. Wang
 \cite{LiWang}, X. Wang \cite{W}, ect..  When the Ricci curvature is super Einstein: $R_{i\jbar}\ge -(n+1) g_{i\jbar}$, Munteanu \cite {M} proves  that 
 $\lambda_1\le n^2$. For the lower bound estimate of $\lambda_1$, Li and Tran \cite{LT} and Li and Wang \cite{LiWang} consider a 
 smoothly bounded pseudoconvex domain in $\CC^n$ with  defining function $r\in C^4(\Dbar)$ such that $u=:-\log (-r)$ is strictly plurisubharmonic
 in $D$. When $r$ is plurisubharmonic in $D$, Li and Tran \cite{LT} prove that $\lambda_1=n^2 $. 
 When $g[u]$ is {\it super asymptotic  Einstein} and $\det H(r)\ge 0$ on $\d D$,
 Li and Wang
 \cite{LiWang} prove $\lambda_1=n^2$.  We will show that $\det H(r)\ge 0$ on $\d D$  when $D$ is super-pseudoconvex. 
 
 The first result of the paper is  the following theorems.

\begin{theorem}  Let $D$ be a smoothly bounded strictly pseudoconvex domain in $\CC^n$.  Let $\tilde{\rho}\in C^4(\Dbar)$ be
a defining function for $D$ such that $\tilde{u}=-\log (-\tilde{\rho})$ is strictly plurisubharmonic.  If the K\"ahler metric $g[\tilde{u}]$ induced
by $\tilde{u}$ is the {\it super asymptotic Einstein}, then  the following two statements hold:

(i) $\tilde{\rho}$ is strictly plurisubharmonic on $\Dbar$ if and only if $D$ is strictly super-pseudoconvex. In particular if
 $\tilde{\rho}=\rho(z)$ is the solution of (1.5) then $\rho$ is strictly plurisubharmonic in $\Dbar$ when $D$ is strictly super-pseudoconvex;
 \smallskip
 
 (ii) If $D$ is also super-pseudoconvex then $\lambda_1(\Delta_{g[\tilde{u}]})=n^2$, where $\Delta_g=-4\sum_{i,j=1}^n g^{i\jbar}{\d^2\over \d z_i \d \zbar_j}$.
\end{theorem}

It is interesting to bridge the relation between convex and super-pseudoconvex. The second result of the paper is:

\begin{theorem} Let $D$ be a smoothly bounded  domain in $\CC^n$. Then

(i) When $n=1$,  $D$ is strictly super-pseudoconvex (super-pseudoconvex) if and only if $D$ is strictly convex (convex);
\smallskip

(ii) When $n>1$,  if $D$ is convex and if there is a strictly plurisubharmonic defining function $ r \in C^4(\Dbar)$ such that
$$
n-1 + { |\d r|^2\over n } a^{k\lbar}[r]\Big[\TD r_{k\lbar}-a^{i\qbar}[r] r^{p\jbar}  r_{i\jbar k} r_{p\qbar \lbar} - (\TD r_k) (\TD r_{\lbar}) \Big]  
 -2\Re  r^k \TD r_k > 0,\leqno(1.16)
$$
then $D$ is strictly super-pseudoconvex;
\smallskip

(iii) Convexity and Super-pseudoconvexity can not contain each other.
\end{theorem}

\medskip

{\bf Acknowledgement:} The author would like to thank Professor C. Fefferman and Xiaodong Wang for some useful conversations he has  had with them.
The author is greatly appreciate  and thank Professor R. Graham who pointed out that there is a mistake in computation at (3.21): $r_{nn}^b=r_{nn}+b r_n^2$ at $z_0$ in the
the previous version of the paper (it should be $r_{nn}^b=r_{nn}+2b r_n^2$),
as well as his valuable suggestions for the current revision. 

\medskip

The paper is organized as follows: Section 2, we give an approximation formula. Theorem 1.3 will be proved in Section 3;
 Part (i) and Part (ii) of  Theorem 1.4 will be proved in Section 4. Finally, in Section 5, we provide two examples which
show that strictly convex and super-pseudoconvex can not contain each other when $n>1$. Which proves Part (iii) of Theorem 1.4.

\section{An approximation formula}

Let $D$ be a bounded domain in $\CC^n$ with smooth boundary.
Let $r\in C^2(\Dbar)$ be a real-valued, negative defining  function for $D$. Then the Fefferman  operator \cite{F1, CY} acting on $r$ is defined by
$$
J(r)=-\det\left[\matrix{ r & \dbar r\cr (\dbar r)^* & H(r)\cr}\right],\leqno(2.1)
$$
where $\dbar r=({\d r\over \d \zbar_1},\cdots, {\d r\over \d \zbar_n})=(r_{\1bar},\cdots, r_{\nbar})$ is a row vector in $\CC^n$ and $(\dbar r)^*$ is
its adjoint vector, which is column vector in $\CC^n$ and $H(r)=[{\d^2 r\over \d z_i \d \zbar_j}]$ is the $n\times n$ complex Hessian matrix of $r$.

If  $H(r)=[r_{i\jbar}]$ is invertible, in particular it is positive definite, then we use the notation $[r^{i\jbar}]^t=:H(r)^{-1}$ and
$$
|\d r|_r^2=\sum_{i,j=1}^n r^{i\jbar} r_i r_{\jbar}.\leqno(2.2)
$$
It is easy to verify that
$$
J(r)=\det H(r) (-r+|\d r|_r^2).\leqno(2.3)
$$
In fact, since
\begin{eqnarray*}
(2.4)\qquad\qquad\qquad\qquad\qquad J(r )&=& (-r) \det [H(r)- {(\dbar r)^* (\dbar r) \over r}]\\
&=&(-r) \det H(r) (1-{|\d r|_r^2 \over r})\\
&=&\det H(r) (-r+|\d r|_r^2).\qquad\qquad\qquad\qquad \qquad\qquad\qquad\qquad 
\end{eqnarray*}

\begin{remark} When $H(r)$ is not positive definite on $\d D$, we can replace $r$ by 
$$
r[a]=:r(z)+{a\over 2} r^2.\leqno(2.5)
$$
 Then $r[a]$ is positive definite with a large $a$ and
$$
J(r)={1\over (1+a r)^n}\det H(r[a]) (-r +(1+2a\, r) |\d r|_{r[a]}).\leqno(2.6)
$$
\end{remark}
\medskip

From now on, we will always assume that $r(z)\in C^\infty(\Dbar)$ be a negative defining function for $D$ such that
$$
\ell(r)=-\log(-r)\leqno(2.7)
$$
is strictly plurisubharmonic in $D$. It is known from \cite{CY, Li1, Li2, Li3} that the following identity holds:
$$
\det H(\ell(r))=J(r) e^{(n+1)\ell(r)}.\leqno(2.8)
$$
This implies that

(i) $u=:\ell (r)$ is strictly plurisubharmonic on $D$ if and only if $J(r)>0$ on $D$;

(ii)  $J(r)=1$ if and only if  $\det H(u)=e^{(n+1) u}$ with $u=:\ell(r)$.
\medskip

C. Fefferman \cite{F1} gave a formula to approximate the potential function $\rho$ (for equation (1.5)).  He proved that $J(r\, J(r)^{-1/(n+1)})=1+O(r)$ near $\d D$. Higher order approximation
can be iterated through the previous steps. Based on the Fefferman's idea, the iteration formula of the approximation was given in  more detail by R. Graham in \cite{G}. The
author \cite{Li1} gave another modification.  For convenience  of  readers and further argument for the current paper, we will state and prove a second order approximation formula here.

\begin{theorem} \sl Let $D$ be a smoothly bounded pseudoconvex domain in $\CC^n$. 
Let $r(z)$ be a smooth negative defining function for $D$ such that $\ell(r)$  is strictly
plurisubharmonic in $D$. Let
$$
\rho_{1}(z) =r(z) J(r)^{-1/(n+1)} e^{-B(z)}\leqno(2.9)
$$
with
$$
B(z)=B[r](z)={ \tr(H(\ell( r ))^{-1} H(\log J(r)) \over 2n (n+1)}.\leqno(2.10)
$$
Then
$$
J(\rho_{1})(z)=1+O(r^2).\leqno(2.11)
$$
Moreover, if $J(r)=1+O(r^2)$ then $\rho_1=r+O(r^3)$ and $J(\rho_1)=1+O(r^3)$.
\end{theorem}

\proof  Since
$$
H(\ell (r))={1\over (-r)(1+ar)} [ H(r_a)+{1+2a\, r \over (-r)} (\dbar r)^*(\dbar r)]\leqno(2.12)
$$
by choosing $a\ge 0$ so that $r[a]$ is strictly plurisubharmonic.  Therefore, we can write
$$
B(z)=(-r) B_0(z),\leqno(2.13)
$$
with $B_0\in C^\infty(\Dbar)$. Since
$$
H(B)=(-r) H(B_0)-B_0 (H(r)+{(\dbar r)^* \dbar r\over -r}) +B_0{(\dbar r)^* \dbar r\over -r}-(\dbar r)^* (\dbar B_0)-(\dbar B)^* (\dbar r).\leqno(2.14)
$$
By complex rotation, one may assume that ${\d r\over \d z_j}(z_0)=0$ for $1\le j\le n-1$ and $H(r)(z_0)$ is diagonal, it is easy to verify that
$$
\tr(H(\ell(r)^{-1}H(B))
=-n B(z)+(-r)B_0+O(r^2)=-(n-1)B +O(r^2).\leqno(2.15)
$$
Since
\begin{eqnarray*}
J(\rho_1)(z)  e^{(n+1)\ell(\rho_1)}&=&\det H(\ell(\rho_1)) \\
&=&\det \Big( H(\ell(r)) +{1\over n+1} H(\log J) +H(B)\Big)\\
&=&\det H(\ell(r)) \det \Big( I_n+  H(\ell(r))^{-1}[{1\over n+1} H(\log J) +H(B)]\Big)\\
&=&J(r)e^{(n+1)\ell(r)} \det \Big( I_n+  H(\ell(r))^{-1}[{1\over n+1} H(\log J) +H(B)]\Big)\\
\end{eqnarray*}
Notice that $\exp((n+1)\ell(\rho_1))=\exp((n+1)B) J(r) \exp((n+1)\ell(r))$, we have
\begin{eqnarray*}
J(\rho_1)(z)  
&=& e^{-(n+1) B} \det \Big( I_n+  H(\ell(r))^{-1}[{1\over n+1} H(\log J) +H(B)]\Big)\\
&=&e^{-(n+1) B}[1+\tr [H(\ell(r))^{-1}[{1\over n+1} H(\log J) +H(B)] +O(r^2)\\
&=&e^{-(n+1) B}[1+2n B +\tr( H(\ell(r))^{-1} H(B)] +O(r^2)\\
&=&e^{-(n+1)B}[1+2n B -(n-1) B+O(r^2)]+O(r^2)\\
&=&1+{(n+1)^2\over 2} B^2+O(r^2)\\
&=&1+O(r^2).
\end{eqnarray*}

When $J(r)=1+A r^2$ with $A$ is smooth on $\Dbar$, it is easy to prove $B=B_1 r^2$ with $B_1$ smooth in $\Dbar$ near $\d D$. It is also easy to verify that
$\rho_1[r]=r+O(r^3)$ and $J(\rho_1[r])=1+O(r^3)$.
This proves Theorem 2.1.\epf

\begin{proposition} Let $D$ be a smoothly bounded strictly pseudoconvex domain
in $\CC^n$.  Let $u$ be the plurisubharmonic solution of (1.3) and
$\rho(z)=-e^{-u}$. Then for any smooth defining function $r$ of $D$ with $\ell(r)$
being strictly plurisubharmonic in $D$, we have
$$
(2.16)\quad \det H(\rho)=J(r)^{-n\over n+1} \det \Big(  H(r) -{  [\d_i r \d_{\jbar} \log J+ \d_i \log J(r)\,  \d_{\jbar} r]\over n+1}
 -[\d_i r \d_{\jbar} B(z)+\d_i B \d_{\jbar} r] \Big)
 $$
on $\d D$,  where $B(z)=B[r](z)$ is given by (2.10).
 \end{proposition}
 
 \proof  Let
 $$
 \rho_1(z)=\rho_1[r]=:r(z) J(r)^{-1/(n+1)} e^{-B}.\leqno(2.17)
 $$
 Theorem 2.1 implies that $\rho(z)=\rho_1(z)+O(r(z)^3)$. A simple calculation shows  that
$$
\det H(\rho)=\det H(\rho_1),\quad z\in \d D.\leqno(2.18)
$$
By (2.13) ($B=(-r) B_0$), one can easily see that
$$
\rho_1(z)=r(z) J(r)^{-1/(n+1)} -r(z) J(r)^{-1/(n+1)} B(z)+O(r(z)^3)\leqno(2.19)
$$
and
$$
\det H(\rho_1)=\det H\Big(r(z) J(r)^{-1/(n+1)} -r(z) J(r)^{-1/(n+1)} B(z)\Big),\quad z\in \d D.\leqno(2.20)
$$
For any $z\in \d D$, by (2.20), one has 
\begin{eqnarray*}
\lefteqn{(2.21)\quad \det H(\rho_1)(z)}\\
&\quad =& \det \Big( H(r J(r)^{-1/(n+1)} ) -J(r)^{-1/(n+1)} [\d_i r \d_{\jbar} B+\d_i B \d_{\jbar} r] \Big)\\
&\quad =& \det \Big( J(r)^{-1\over (n+1)}  H(r) -{ J^{-(n+2)\over (n+1)}\over n+1} [\d_i r \d_{\jbar} J+ \d_i J(r) \d_{\jbar} r]
 -J(r)^{-1\over (n+1)} [\d_i r \d_{\jbar} B+\d_i B \d_{\jbar} r] \Big)\\
 &\quad =&J(r)^{-n\over n+1} \det \Big(  H(r) -{ 1\over n+1} [\d_i r \d_{\jbar} \log J+ \d_i \log J(r)\,  \d_{\jbar} r]
 -[\d_i r \d_{\jbar} B+\d_i B \d_{\jbar} r] \Big).
\end{eqnarray*}
This proves Proposition 2.2.\epf
\medskip

Let $u^{D_j}$ be the potential functions for the K\"ahler-Einstein metric for $D_j$ and let
$$
\rho^{D_j}(z)=-e^{-u^{D_j}(z)},\quad j=1, 2.\leqno(2.22)
$$
\begin{proposition} Let $\phi: D_1\to D_2$ be a smooth biholomorphic mapping. Then
$$
\rho^{D_1}(z)=\rho^{D_2}(\phi(z)) |\det\phi'(z)|^{-2/(n+1)}\leqno(2.23)
$$
In particular, if $\det\phi'(z)$ is constant $c$ then
$$
\det H(\rho^{D_1})(z)=|c|^{2/(n+1)} \det H(\rho^{D_2})(\phi(z)). \leqno(2.24)
$$
\end{proposition}

\proof Since $\phi: D_1\to D_2$ is biholomorphic, one has that if $u^{D_j}$ is the unique plurisubharmonic solutions for
the Monge-Amp\`ere equation:
$$
\cases{\det H(u)=e^{(n+1) u}, \quad z\in D_j\cr \qquad \qquad u=\infty,\quad z\in \d D_j\cr}\leqno(2.25)
$$
Then 
$$
u^{D_1}(z)=u^{D_2}(\phi(z))+{1\over n+1}\log |\det \phi'(z)|^2,\quad z\in D_1\leqno(2.26)
$$
and
$$
\rho^{D_1}(z)=\rho^{D_2}(\phi(z))|\det\phi'(z)|^{-2/(n+1)}.\leqno(2.27)
$$
In particular, when $\det \phi'(z)=c$, one has
$$
\det H(\rho^{D_1})(z)=|c|^{-2n/(n+1)} \det H(\rho^{D_2})(\phi(z)) |c|^2=|c|^{2/(n+1)} \det H(\rho^{D_2})(\phi(z)) 
$$
and the proof of Proposition 2.3 is complete.\epf
\medskip

We also need the following holomorphic change of variables formula.

\begin{lemma} For $z_0\in \d D$, if  $z=\phi(w): B(0,\delta_0) \to B(z_0, 1)$ be a  one-to-one holomorphic map with $\phi(0)=z_0$ and
$r(z)=\tilde{r}(w)$, then 
%$$
%\rho_1(\phi(w) )=|\det \phi'(z)|^{2/(n+1)} {r(\phi)\over J(r\circ\phi)^{1/(n+1)}}e^{-B(r\circ\phi)}.\leqno(2.27)
%$$
%If $z=\psi(w)$, then 
$$
\rho_1(\phi(w))=|\det \phi'(w)|^{2/(n+1)} {\tilde{r}(w)\over J(\tilde{r}(w))^{1/(n+1)}}e^{-B(\tilde{r}(w))}.\leqno(2.28)
$$
Moreover,  if  $|\det \phi'(z)|^2$ is a constant
on $B(0,\delta_0)$ for some $\delta_0>0$
$$
\det H(\rho_1)(z_0) |\det \phi'(0)|^{2\over n+1}= \det H\Big( {\tilde{r}\over J(\tilde{r})^{1/(n+1)}}e^{-B(\tilde{r})}\Big)(0).\leqno(2.29)
$$
\end{lemma}

\proof Since $|\det\phi'(z)|^2$ is constant, by the definitions for $B[r]$ and $J(r)$ from Theorem 2.1, one can easily prove (2.27) and (2.29), 
and  the proposition is proved.\epf

\section{Proof of Theorem 1.3}

Let $D$ be a smoothly bounded strictly pseudoconvex domain in $\CC^n$. 
Let $r\in C^\infty(\Dbar)$ be any strictly plurisubharmonic defining function for $D$. Let
 $$
 \rho_1(z)= r(z) J(r)^{-1/(n+1)} \exp (-B(z)) \leqno(3.1)
 $$
 where
 $$
 B(z)={ \tr(H(\ell(r))^{-1} H(\log J(r)) \over 2n (n+1)},\leqno(3.2)
 $$
According to Theorem 2.1, one has
$$
J(\rho_1)=1+O(r(z)^2).\leqno(3.3)
$$
Let $\rho=\rho^D$ be the solution of (1.5) such that $\ell(\rho)$ is strictly plurisubharmonic in $D$. Then
$$
\det H(\rho)(z)=\det H(\rho_1)(z)\quad\hbox{on }\d D.\leqno(3.4)
$$

\medskip

By Proposition 2.2 and
\begin{eqnarray*}
(3.5)\qquad \qquad B(z)&=&{(-r)\over 2n (n+1)} \tr [(H(r)+{r_i r_{\jbar}\over -r})^{-1} H(\log J(r) ] (z)\\
&=& {(-r)\over 2n (n+1)}\sum_{j,k=1}^n (r^{i\jbar} -{r^i r^{\jbar} \over -r+|\d r|_r^2}) {\d^2 \log J(r)\over \d z_i\d \zbar_j}\qquad\qquad\qquad\qquad\\
&=&-B^0(z) r,
\end{eqnarray*}
where
\def\TD{\tilde{\Delta}}
$$
B^0(z)={1\over 2n (n+1)} \sum_{j,k=1}^n a^{i\jbar}[r] {\d^2 \log J(r)\over \d z_i\d \zbar_j}={1\over 2n(n+1)}{\TD}_r \log J(r).\leqno(3.6)
$$
Thus for $z_0\in \d D$, one has 
$$
\d_j B (z_0)=- B^0(z_0) \d_j r (z_0), \quad \d_{\jbar} B (z_0)=- B^0(z_0) \d_{\jbar} r (z_0), \quad\hbox{ for } 1\le j\le n.\leqno(3.7)
$$
Let
$$
R=\sum_{j=1}^n r^{j} {\d\over \d z_j},\quad \Rbar=\sum_{j=1}^n r^{\jbar}{\d\over \d \zbar_j},\quad r^i=r^{i\jbar} r_{\jbar}, \ \ r^{\jbar}=r^{i\jbar} r_i.\leqno(3.8)
$$
and
$$
\Big| \tilde{\nabla_r} f \Big |^2=:\sum_{i,j=1}^n (r^{i\jbar}-{r^i r^{\jbar}\over -r +|\d r|_r^2})\d_i f\d_{\jbar} f=\sum_{i,j=1}^n r^{i\jbar} \d_i f\d_{\jbar} f-{|R f|^2\over -r +|\d r|^2_r }.\leqno(3.9)
$$
Then it is easy to see that
$$
|\tilde{\nabla}_r r|^2=0\quad \hbox{ on }\  \d D.\leqno(3.10)
$$
Therefore, by (2.21) and Lemma 3.1 in \cite{Li1}, at $z=z_0\in \d D$, one has
\begin{eqnarray*}
\lefteqn{(3.11)\qquad \det H(\rho)(z^0)\,J(r)^{n/(n+1)}(z^0)}\\
&\qquad  =&\det H(r) \Big(\Big| 1- r^{i\jbar} (\d_i r\, ({\d_{\jbar}\log J(r) \over n+1}-B^0 \d_{\jbar} r)\Big|^2\\
&&- |\d r|_r^2 \sum_{i,j=1}^n r^{i\jbar}
({\d_{i}\log J(r) \over n+1}-B^0 \d_i r) ({\d_{\jbar}\log J(r) \over n+1}-B^0 \d_{\jbar} r)\Big)\\
&\qquad =&\det H(r) \Big(\Big| |1- { \Rbar \log J(r) \over n+1}+B^0 |\d r|^2_r\Big|^2\\
&&- |\d r|_r^2 \sum_{i,j=1}^n r^{i\jbar}
{\d_{i}\log J(r) \d_{\jbar } \log J(r)\over (n+1)^2} +|\d r|_r^2 2\Re  B^0 {\Rbar \log J(r)\over n+1} - |\d r|_r^4 |B^0|^2  \Big)\\
&\qquad =&\det H(r) \Big(1+2 B^0 |\d r|^2 -2\Re {\Rbar \log J(r)\over n+1} -{|\d r|_r ^2\over (n+1)^2} |\tilde{\nabla}_r \log J(r)|^2  \Big)\\
&\qquad =&\det H(r) \Big(1+{|\d r|^2\over n(n+1)} \TD \log J(r) -2\Re {\Rbar \log J(r)\over n+1} -{|\d r|_r ^2\over (n+1)^2} |\tilde{\nabla}_r \log J(r)|^2  \Big)\qquad\\
&\qquad >& 0
\end{eqnarray*}
since $D$ is strictly super-pseudoconvex, there is a strictly plurisubharmonic function $r\in C^4(\overline{D})$ such that the above inequality holds on $\d D$.
If $\tilde\rho$ is smooth defining function for $D$ such that the K\"ahler metric induced by $\tilde{u}=-\log(-\tilde{\rho})$ is super asymptotic Enistein, then
$\det H(\tilde{\rho})=\det H(\rho) >0$ on $\d D$ by (3.11). By Lemma 2 in \cite{LiWang}, one has that
 $\det H(\tilde{\rho})$ attains its minimum over $\overline{D}$ at some ponit in $\d D$.
  Therefore, $\det H(\tilde{\rho})>0$ on $\Dbar$ and the proof of  Part (i) of Theorem 1.3 is complete.
Part (ii) of Theorem 1.3 is a corollary of Part (i)  and the result in \cite{LT} and \cite{LiWang}. Therefore, the proof of Theorem 1.3 is complete. \epf

\section{Super-pseudoconvex domains}

In this section we will study more on the super-pseudoconvex domain in $\CC^n$ comparing with convex domains.
Since 
$$
\log J(r)=\log \det H(r)+\log (-r+|\d r|_r^2),\leqno(4.1)
$$
\begin{eqnarray*}
(4.2)\qquad {\d (-r +|\d r|^2_r)\over \d z_k}&=&-r_k+ \d_k(r^{i\jbar})  r_ i r_{\jbar} +r^{i\jbar} r_{ik} r_{\jbar}+r^{i\jbar} r_i r_{k\jbar} \qquad\qquad\qquad\qquad\\
&\qquad =& -r ^{i \qbar} r^{p\jbar} r_{p\qbar k}   r_ i r_{\jbar} +r^{i\jbar} r_{ik} r_{\jbar} \\
&\qquad =& -r ^{ \qbar} r^{p} r_{p\qbar k}   +r^{i} r_{ik} 
\end{eqnarray*}
and
$$
{\d \log J\over \d z_k}={\d \log \det H(r)+\log (-r+|\d r|^2_r)\over \d z_k}=(r^{i\jbar} -{r^i r^{\jbar}\over -r+|\d r|_r^2} ) r_{i\jbar k} +{r^i r_{ik}\over -r+|\d r|_r^2},\leqno(4.3)
$$
we have
$$
R \log J(r)(z_0)= r^k {\TD} r_k + {r^i r^k\over |\d r|_r^2} r_{ik}.\leqno(4.4)
$$
Thus,
$$
 \det H(\rho)(z^0)\,J(r)^{n/(n+1)}(z^0)
=\det H(r) \Big(1-{2 \Re r^k r^i r_{ik}\over (n+1)|\d r|^2}+ \tilde{E}(r)\Big),\leqno(4.5)
$$
where
$$
\tilde{E}(r)=: {  |\d r|^2 \over n(n+1)}\Big[ {\TD} \log J(r) -{n |{\TN}\log J(r)|^2\over (n+1)}  -2n \Re ({r^k  \TD r_k\over |\d r|_r^2}  )\Big]. \leqno(4.6)
$$

\begin{proposition} Let $D$ be a smoothly bounded domain in the complex plane $\CC$. Then  $D$ is (strictly) super-pseudoconvex 
if and only if $D$ is (strictly) convex.
\end{proposition}

\proof  Let $r$ be any smooth strictly subharmonic defining function on $D\subset \CC$. By (4.5) and (4.6),  we have
$a^{1\1bar}[r]=0$ and $\tilde{E}(r)=0$ on $\d D$. Therefore,  $D$ is strictly super-pseudoconvex if and only if
$$
 S_r(z)=:\det H(r) \Big(1-{2\over n+1}\Re {r^k r^i r_{ik}\over |\d r|_r^2}\Big) >0 \leqno(4.7)
$$
on $\d D$. For ant $z_0\in \d D$, by rotation, we may assume that $ r_n(z_0)>0$. Thus
$$
S_r(z_0)
=r_{1\1bar} -\Re r_{11}(z_0)\leqno(4.8)
$$
is positive for all $z_0 \in \d D$ if and only if $\d D$ is strictly convex; and is non-negative for all $z_0\in \d D$  if and only if $\d D$ is convex, respectively. Therefore,
the proof of the proposition is complete. \epf
\medskip

\noindent{\bf Next we estimate $\tilde{E}(r)$. }
\medskip
\begin{proposition} With the notation above, for $z\in \d D$, we have
$$
\tilde{E}(r)\ge  { |\d r|^2 a^{k\lbar}[r]\over n (n+1)} \Big[\TD r_{k\lbar}-a^{i\qbar}[r] r^{p\jbar}  r_{i\jbar k} r_{p\qbar \lbar} - (\TD r_k) (\TD r_{\lbar})  
- n
 {r^i r_{ik} r^{\jbar} r_{\jbar \lbar}\over |\d r|_r^4}\Big] 
 -{2\Re  r^k \TD r_k  \over (n+1) }.\leqno(4.9)
$$
and
$$
\tilde{E}(r)\le   {|\d r|^2 a^{k\lbar} \over n (n+1)}\Big[\TD r_{k\lbar}+a^{i\qbar}[r] r^p r^{\jbar}r_{i\jbar k}r_{p\qbar\lbar}+2 a^{i\qbar}[r] \, {r_{i k} r_{\qbar \lbar} \over |\d r|^2} \Big]-{2\Re  r^k \TD r_k  \over (n+1) }.\leqno(4.10)
$$
\end{proposition}
\proof Notice
$$
(r^i)_{\lbar}=(r^{i\qbar} r_{\qbar})_{\lbar}=r_{\qbar} (r^{i\qbar})_{\lbar}+r^{i\qbar} r_{\qbar \lbar}=-r^{i \tbar} r^{s\qbar} r_{s \tbar \lbar} r_{\qbar}+r^{i \qbar} r_{\qbar \lbar}
=-r^{i\tbar} r^s r_{s\tbar \lbar}+r^{i\qbar} r_{\qbar\lbar}
$$
and
$$
(r^{\jbar})_{\lbar}=(r^{p \jbar} r_p)_{\lbar}=-r^{\qbar} r^{i\jbar} r_{i \qbar \lbar}+\delta_{j\ell}.
$$
By (4.3) and (4.2), for $z\in \d D$, one has
\begin{eqnarray*}
{\d^2 \log J(r)\over \d z_k \d \zbar_{\ell}}&=&(r^{i\jbar} -{r^i r^{\jbar}\over |\d r|_r^2} ) r_{i\jbar k\lbar} +r_{i\jbar k}{\d \over \d \zbar_\ell} (r^{i\jbar}-{r^i r^{\jbar}\over -r +|\d r|^2_r})
+{\d \over \d \zbar_\ell}{r^i r_{ik}\over (-r+|\d r|_r^2)}\\
&=&\TD r_{ k\lbar} -r_{i\jbar k} r^{i\qbar} r^{p\jbar} r_{p\qbar\lbar}\\
&&+{1\over (|\d r|_r^2)^2} (r_{i\jbar k}r^i r^{\jbar}
-r^i r_{ik})( {\d (-r +|\d r|_r^2) \over \d \zbar_\ell})\\
&&-{r_{i\jbar k}\over |\d r|^2_r}  (r^i (r^{\jbar})_{\lbar} +r^{\jbar} (r^i)_{\lbar})+ {1\over |\d r|^2} (r^i r_{ik\lbar}+ r_{ik} (r^i)_{\lbar})\\
&=&\TD r_{ k\lbar} -r_{i\jbar k} r^{i\qbar} r^{p\jbar} r_{p\qbar\lbar}\\
&&+{1\over (|\d r|_r^2)^2} (r_{i\jbar k}r^i r^{\jbar}
-r^i r_{ik})(  -r ^{ \qbar} r^{p} r_{p\qbar \lbar}   + r^{\qbar} r_{\qbar \lbar} )\\
&&-{r_{i\jbar k}\over |\d r|^2_r}  \Big( r^{\jbar}( -r^{i\tbar} r^s r_{s\tbar \lbar}+r^{i\qbar} r_{\qbar\lbar}) +r^i (-r^{\qbar} r^{p\jbar} r_{p \qbar \lbar}+\delta_{j\ell})\Big)\\
&&+ {1\over |\d r|^2} \Big(r^i r_{ik\lbar}+ r_{ik} (-r^{i\tbar} r^s r_{s\tbar \lbar}+r^{i\qbar} r_{\qbar\lbar})\Big)\\
&=&\TD r_{k\lbar}-r^{i\qbar} r^{p\jbar} r_{i\jbar k} r_{p\qbar \lbar} -{1\over (|\d r|_r^2)^2} (r_{i\jbar k}r^i r^{\jbar}
-r^i r_{ik})(  r ^{ \qbar} r^{p} r_{p\qbar \lbar}   - r^{\qbar} r_{\qbar \lbar} )\\
&&+{1 \over |\d r|^2_r} (r^p r^{\jbar} r^{i\qbar} +r^i r^{\qbar} r^{p\jbar})  r_{p\qbar \lbar}r_{i\jbar k}-{1\over |\d r|_r^2} r^{\jbar} r^{i\qbar} r_{\qbar\lbar}r_{i\jbar k} -{r_{i\lbar k}\over |\d r|_r^2} r^i\\
&&+ {1\over |\d r|^2} \Big(r^i r_{ik\lbar}-r^{i\tbar} r^s r_{s\tbar \lbar} r_{ik} +r^{i\qbar} r_{\qbar\lbar} r_{ik} )\\
&=&\TD r_{k\lbar}-r^{i\qbar} r^{p\jbar} r_{i\jbar k} r_{p\qbar \lbar} -{r^i r^{\jbar} r^p r^{\qbar} \over |\d r|_r^4} r_{i\jbar k} r_{p\qbar \lbar} +{1\over |\d r|_r^4}(r^i r^{\jbar} r_{i\jbar k} r^{\qbar} r_{\qbar\lbar}
+r^p r^{\qbar} r_{p\qbar \lbar} r^i r_{ik})
 \\
&&+{1 \over |\d r|^2_r} (r^p r^{\jbar} r^{i\qbar} +r^i r^{\qbar} r^{p\jbar})  r_{p\qbar \lbar}r_{i\jbar k}\\
&&- {1\over |\d r_r^2} \Big( r^i r^{p\jbar} r_{pk} r_{i\jbar \lbar} +r^{\jbar} r^{i\qbar} r_{\qbar \lbar} r_{i\jbar k} ) +{1\over |\d r|_r^2} (r^{i\qbar} -{r^i r^{\qbar}\over |\d r|_r^2}) r_{\qbar\lbar} r_{ik} \\
&=&\TD r_{k\lbar}-(r^{i\qbar}-{r^i r^{\qbar}\over |\d r|_r^2}) (r^{p\jbar} -{r^p r^{\jbar}\over |\d r|_r^2}) r_{i\jbar k} r_{p\qbar \lbar}   \\
%&&+{1 \over |\d r|^2_r} r^i r^{\qbar} (r^{p\jbar} -{r^p r^{\jbar}\over |\d r|_r^2})  r_{p\qbar \lbar}r_{i\jbar k}\\
&&- {1\over |\d r|_r^2} \Big( r^i ( r^{p\jbar} -{r^p r^{\jbar}\over |\d r|_r^2} ) r_{pk} r_{i\jbar \lbar} +r^{\jbar} (r^{i\qbar} -{r^i r^{\qbar}\over |\d r|_r^2}) r_{\qbar \lbar} r_{i\jbar k} \Big) +{1\over |\d r|_r^2} (r^{i\qbar} -{r^i r^{\qbar}\over |\d r|_r^2}) r_{\qbar\lbar} r_{ik}.
\end{eqnarray*}
Then for $z\in \d D$, we have
\begin{eqnarray*}
\TD \log J(r)(z)&\ge & a^{k\lbar}[r]\TD r_{k\lbar}-a^{k\lbar}[r] a^{i\qbar}[r] a^{p\jbar}[r]  r_{i\jbar k} r_{p\qbar \lbar}   \\
&&-a^{k\lbar}[r]  {a^{i\qbar}[r] \over |\d r|^2_r} \Big( r^{\jbar} r_{i\jbar k} r^p r_{p\qbar\lbar} +r_{ki} r_{\qbar \lbar}\Big) +{1\over |\d r|_r^2} a^{k\lbar}[r] a^{i\qbar}[r]  r_{\qbar\lbar} r_{ik} \\
&=&a^{k\lbar}\TD r_{k\lbar}-a^{k \lbar}[r] a^{i\qbar}[r] \, r^{p\jbar}  r_{i\jbar k} r_{p\qbar \lbar}   \\
\end{eqnarray*}
and
\begin{eqnarray*}
\TD \log J(r)(z)&\le &a^{k\lbar}\TD r_{k\lbar}+2a^{k \lbar}[r] a^{i\qbar}[r] \, {r_{i k} r_{\qbar \lbar} \over |\d r|^2} +a^{k\lbar}[r] a^{i\qbar}[r] r^p r^{\jbar} r_{i\jbar k} r_{p\qbar\lbar} \\
\end{eqnarray*}
Moreover,
\begin{eqnarray*}
|\TN \log J(r)|^2&=& a^{k\lbar}[r]\, (\TD r_k +{r^i r_{ik}\over |\d r|_r^2})(\TD r_{\lbar}+{r^{\jbar} r_{\jbar \lbar}\over |\d r|_r^2})\\
&=& a^{k\lbar}[r]\, \Big[(\TD r_k) (\TD r_{\lbar}) +(\TD r_k)({r^{\jbar} r_{\jbar \lbar}\over |\d r|_r^2})
 +{r^i r_{ik}\over |\d r|_r^2} \TD r_{\lbar}+{r^i r_{ik}\over |\d r|_r^2} {r^{\jbar} r_{\jbar \lbar}\over |\d r|_r^2}\Big]\\
 &\le& a^{k\lbar}[r]\, \Big[{n+1\over n} (\TD r_k) (\TD r_{\lbar}) 
+(n+1){r^i r_{ik}\over |\d r|_r^2} {r^{\jbar} r_{\jbar \lbar}\over |\d r|_r^2}\Big].
\end{eqnarray*}
Therefore,
\begin{eqnarray*}
\lefteqn{\TD \log J(r)-{n\over n+1} |\TN \log J|^2}\\
&\ge &a^{k\lbar}[r]\Big(\TD r_{k\lbar}-a^{i\qbar}[r] r^{p\jbar}  r_{i\jbar k} r_{p\qbar \lbar} \Big) 
-a^{k\lbar}[r]\Big ( (\TD r_k) (\TD r_{\lbar}) 
+n {r^i r_{ik}\over |\d r|_r^2} {r^{\jbar} r_{\jbar \lbar}\over |\d r|_r^2}\Big).
\end{eqnarray*}
Therefore,
\begin{eqnarray*}
\tilde{E}(r)&\ge& { |\d r|^2 a^{k\lbar}[r]\over n (n+1)}\Big(\TD r_{k\lbar}-a^{i\qbar}[r] r^{p\jbar}  r_{i\jbar k} r_{p\qbar \lbar} 
- (\TD r_k) (\TD r_{\lbar}) 
-n {r^i r_{ik}\over |\d r|_r^2} {r^{\jbar} r_{\jbar \lbar}\over |\d r|_r^2}\Big)  
 -{2\Re  r^k \TD r_k  \over (n+1) }.
\end{eqnarray*}
and
$$
\tilde{E}(r)\le {|\d r|^2 a^{k\lbar}[r] \over n (n+1)}\Big[\TD r_{k\lbar}+ a^{i\qbar} r^p r^{\jbar} r_{i\jbar k} r_{p\qbar \lbar}+ 2a^{i\qbar}[r] \, {r_{i k} r_{\qbar \lbar} \over |\d r|^2} \Big]-{2\Re  r^k \TD r_k  \over (n+1) }
$$
Therefore, the proof of the proposition is complete.\epf

\begin{corollary} Let $D$ be smoothly bounded convex domain in $\CC^n$. If there is a strictly plurisubharmonic defining function $r\in C^4(\Dbar)$ 
such that 
$$
{n-1\over n+1}+ { |\d r|^2 a^{k\lbar}[r]\over n (n+1)}\Big(\TD r_{k\lbar}-a^{i\qbar}[r] r^{p\jbar}  r_{i\jbar k} r_{p\qbar \lbar} 
- (\TD r_k) (\TD r_{\lbar}) 
\Big)  
 -{2\Re  r^k \TD r_k  \over (n+1) } >0
   \hbox{ on }\ \d D, \leqno(4.11)
$$
then $D$ is strictly super-pseudoconvex.
\end{corollary}

\proof If $\d D$ is convex then for any strictly plurisubharmonic defining function $r\in C^4(\Dbar)$, we have
$$
{2\over n+1}-{2\over n+1}\Re { r^k r^i r_{ik}\over |\d r|^2} -
 {a^{k\lbar}[r] r^i r_{ik} r^{\jbar} r_{\jbar \lbar}\over (n+1) |\d r|_r^2}\ge 0 \quad\hbox{on }\ \d D.\leqno(4.12)
 $$
 Since 
 $$
\tilde{E}(r) +{1\over n+1} a^{k\lbar}[r] r^i r_{ik} r^{\jbar} r_{\jbar\lbar} =  { |\d r|^2 a^{k\lbar}[r] \over n (n+1)}\Big(\TD r_{k\lbar}-a^{i\qbar}[r] r^{p\jbar}  r_{i\jbar k} r_{p\qbar \lbar} 
- (\TD r_k) (\TD r_{\lbar}) \Big)  
 -{2\Re  r^k \TD r_k  \over (n+1) }
 $$
 and 
 $1-{2\over n+1}={n-1\over n+1}$,  by  (4.5), (4.11) and  (4.12), we have  $\det H(\rho)>0$ on $\d D$. This implies 
 $\rho$ is strictly plurisubharmonic on $\Dbar$ by Lemma 2 in \cite{LiWang}.
 This proves Parts (i) and (ii) in Theorem 1.4.\epf

\section{Examples}

In this section, we will provide two examples which show that strictly convex domain and strictly super-pseudoconvex can not contain each other.
\medskip

For $\delta=4^{-12}$, we let
$$
g(t)=:g_{\delta} (t)=:\cases{ e^{-{\delta\over \delta-t}}, \quad \hbox{ if } t<\delta,\cr 0,\quad\qquad  \hbox{ if  } t \ge \delta.\cr}\leqno(5.1)
$$
Let
$$
r(z)=-2\Re z_2+|z|^2-8 |z_1|^4 g(|z_1|^2),\quad z=(z_1, z_2)\in \CC^2.\leqno(5.2)
$$

\begin{example} Let $D=\{z\in \CC^2: r(z)<0\}$. Then

(i)  $D$ is strictly convex.

(ii)  If  $\rho_D$ the solution
of Fefferman equation (2), then $\rho_D$ is not plurisubharmonic in $D$.  

\end{example}

\proof Since
$$
{\d |z_1|^4 g(|z_1|^2)\over \d x_1}=4|z_1|^2 x_1 g(|z_1|^2)+|z_1|^4 g'(|z_1|^2) 2x_1,
$$
$$
 {\d |z_1|^4 g(|z_1|^2)\over \d y_1}=4|z_1|^2 y_1 g(|z_1|^2)+|z_1|^4 g'(|z_1|^2) 2y_1,
$$
$$
{\d^2 \,  |z_1|^4 g(|z_1|^2)\over \d x_1^2 }=16 |z_1|^2 x_1^2 g'(|z_1|^2) +2|z_1|^4 g'(|z_1|^2)+ 4(|z_1|^2 +2 x_1^2)  g(|z_1|^2) + 4|z_1|^4 g''(|z_1|^2) x_1^2,
$$
$$
 {\d |z_1|^4 g(|z_1|)\over \d y_1^2}=16 |z_1|^2 y_1^2 g'(|z_1|^2) +2|z_1|^4 g'(|z_1|^2)+ 4(|z_1|^2 +2 y_1^2) g(|z_1|^2) + 4|z_1|^4 g''(|z_1|^2) y_1^2,
$$
and
\begin{eqnarray*}
{\d^2 (|z_1|^4 g(|z_1|^2)\over \d x_1 \d y_1}
&=&{\d (4|z_1|^2 x_1 g(|z_1|^2)+|z_1|^4 g'(|z_1|^2) 2x_1)\over \d y_1}\\
&=&8 x_1y_1 g(|z_1|^2)+16|z_1|^2 x_1 y_1 g'(|z_1|^2)+4|z_1|^4x_1 y_1 g''(|z_1|^2)
\end{eqnarray*}
Since
\begin{eqnarray*}
20 t^2 |g'(t)|+12 t g(t)+4t^3 |g''(t)|
&=&4t g(t) [3 +5 {t \delta\over (\delta-t)^2}+{t^2 (\delta^2+2\delta (\delta-t))\over (\delta-t)^4}]\\
&\le & 4t g(t) [{11\delta^4\over (\delta-t)^4}]\\
&\le& 4^7 \delta\\
&\le &4^{-5}
\end{eqnarray*}
This implies
$$
18|z_1|^4|g'(|z_1|^2)|+12|z_1|^2 g(|z_1|^2)+4|z_1|^6 |g''(|z_1|^2)|\le {1/4}
$$
and
$$
\Big|{\d (|z_1|^4 g(|z_1|^2)\over \d x_1^2}\Big|<1/4,\ \Big|{\d (|z_1|^4 g(|z_1|^2)\over \d y_1^2}\Big|<1/4 \quad\hbox{and } 
\Big|{\d (|z_1|^4 g(|z_1|^2)\over \d x_1\d y_1}\Big|<1/2
$$
Then 
$
D^2 r(z)=2I_n+D^2 (|z_1|^4 g(|z_1|^2))$ is positive definite in $\RR^4$. Therefore, $D$ is strictly convex. Moreover, $H(r)(0)=I_2$. We claim that
$$
\det H(\rho_D)(0)<0.
$$
Since, at $z=0$, we have
$$
{\d r\over \d z_2}=-1,\quad r_{kj}(0)=r_{i\jbar k}(0)=0,\quad 1\le i,j,k\le 2
$$
By (4.3). This implies ${\d \log J(r )\over \d z_j}(0)=0$ for all $1\le j\le 2$. By (4.6) and (4.10), we have
$$
r_{1\1bar1\1bar}(0)=-32e^{-1},\quad \tilde{E}(r)(0)={|\d r|^2 \over 6 } r_{1\1bar 1\1bar}(0)=-{32\over 6} e^{-1}
$$
Thus,
$$
\det H(\rho_D) J(r)^{2/3} =1-{2\over 3}-{32\over 6 e}<0.
$$
This completes the proof of the statement in the example.\epf
\bigskip

\begin{example} For $n\ge 2, \alpha =21/20$ and $0<C\le (9-8\alpha)(1+\alpha)/256$, we let $r(z)=|z|^2+2\Re z_n 
+\alpha\Re \sum_{j=1}^n z_j^2 +C\sum_{j=1}^n |z_j|^4$ and let
$$
D=\{z\in \CC^n: r(z)<0\}
$$
Then $D$ is super-pseudoconvex, but $D$ is not convex
\end{example}

\proof At $z=(0,0,\cdots, 0)\in \d D$, we have that ${\d\over \d x_j}, {\d \over \d y_j}$ and ${\d \over \d y_n}$ are tangent vectors to $\d D$ for
$1\le j\le n-1$.  Since
$$
{\d^2 r\over \d y_n^2}=2-2\alpha=-2(\alpha-1)<0.
$$
It is easy to see that  $\d D$ at $z=0$, and so $\d D$  is not convex. However,
$$
H(r)=I_n+4C\hbox{Diag}(|z_1|^2,\cdots, |z_n|^2)
$$
where $\hbox{Diag}(|z_1|^2,\cdots, |z_n|^2)$ is a diagonal matrix with diagonal entries $|z_1|^2,\cdots, |z_n|^2$, respectively. Then
$$
{\d^2 r \over \d z_i \d \zbar_j \d z_k \d \zbar_\ell}(z)=4C \delta_{ij} \delta_{k\ell}\delta_{ik},\quad {\d^3 r\over \d z_k \d \zbar_\ell \d z_j}=4C\delta_{k\ell} \delta_{kj} \zbar_j,
\quad {\d^2 r\over \d z_i \d z_j}=( \alpha +2C  \zbar_j^2) \delta_{ij}.
$$
For each $i$
$$
r^i={r_{\ibar}\over 1+4C|z_i|^2},\quad |\d r|_r^2=r^i r_i=\sum_{i=1}^n {|r_i|^2\over 1+4C|z_i|^2}
$$
and, on $\d D$, we have
$$
\TD=\sum_{i, j=1}^n ({\delta_{ij}\over 1+4C|z_j|^2}-{r_{\ibar} r_j\over (1+4C|z_i|^2)(1+4C|z_j|^2)  |\d r|_r^2}){\d^2 \over \d z_i \d \zbar_j}
$$
Notice that if $z\in D$, we have 
$$
2x_n+(1+\alpha)\sum_{j=1}^n x_j^2+(1-\alpha)\sum_{j=1}^n y_j^2+C\sum (x_j^2+y_j^2)^2<0.
$$
This implies that
$$
2x_n+(1+\alpha) x_n^2<0\iff -{2\over 1+\alpha} <x_n<0.\leqno(5.1)
$$
Thus
$$
2x_n+(1+\alpha)x_n^2>{-1\over 1+\alpha} \quad\hbox{and } \ C|z_k|^4-(\alpha-1)|z_k|^2<{1\over 1+\alpha}.\leqno(5.2)
$$
We claim that 
$$
4C|z_k|^2 \le 1/8 \quad\hbox{if }\ 0<C\le {(9-8\alpha)(1+\alpha)\over 256}, 1<\alpha<9/8.\leqno(5.3)
$$
Otherwise, $4C|z_k|^2\ge 1/8$. Then
$
C|z_k|^4-(\alpha-1)|z_k|^2<{1\over 1+\alpha} $ implies 
$$
 |z_k|^2<{8 \over (1+\alpha)(9-8\alpha)}.
$$
This is a contradiction with $4C|z_k|\ge 1/8$. Therefore, the claim is true. Notice 
$$
a^{k\lbar}[r] r_{\lbar}=0, \hbox{ for all } 1\le k\le n,
$$
we have
\begin{eqnarray*}
(r^{k\lbar}-{r^k r^{\lbar}\over  |\d r|^2})(r^i r_{ik} r^{\jbar} r_{\jbar \lbar})
&=&(r^{k\lbar}-{r^k r^{\lbar}\over |\d r|^2}) r^k r^{\lbar} (\alpha+2C\zbar_k^2)(\alpha+2Cz_{\ell}^2)\\
&=&(r^{k\lbar}-{r^k r^{\lbar}\over |\d r|^2}) r_k r_{\lbar} (\alpha+2C{\zbar_k^2-2\alpha |z_k|^2\over 1+4C|z_k|^2})(\alpha+2C{z_{\ell}^2-2\alpha |z_\ell|^2\over 1+4C|z_\ell|^2})\\
%&=&(r^{k\lbar}-{r^k r^{\lbar}\over |\d r|^2}) r_k r_{\lbar} \alpha^2   r^{k\kbar} (1+{2\over \alpha} C\zbar_k^2)r^{\ell\lbar} (1+{2\over \alpha} Cz_{\ell}^2)\\
&=&(r^{k\lbar}-{r^k r^{\lbar}\over |\d r|^2}) r_k r_{\lbar}\alpha^2  \\
&+&4C\alpha \Re (r^{k\lbar}-{r^k r^{\lbar}\over |\d r|^2}) r_k r_{\lbar}  {\zbar_k^2-2\alpha|z_k|^2\over 1+4C|z_k|^2}) \\
&+&4C^2 (r^{k\lbar}-{r^k r^{\lbar}\over |\d r|^2}) r_k r_{\lbar} {(\zbar_k^2-2\alpha |z_k|^2) (z_{\ell}^2-2\alpha |z_{\ell}|^2)\over (1+4C|z_k|^2)(1+4C|z_{\ell}|^2)}\\
%&=&4C^2 (r^{k\lbar}-{r^k r^{\lbar}\over |\d r|^2}) r_k r_{\lbar} r^{k\kbar} (\zbar_k^2-2\alpha |z_k|^2]  r^{\ell\lbar} (z_{\ell}^2
%-2\alpha |z_{\ell}|^2]\\
&\le&{4C^2 (2\alpha+1)^2 |z_k|^4\over (1+4C|z_k|^2)^2} r^{k\kbar} |r_k|^2\\
&\le& {(2\alpha+1)^2 \over 256}|\d r|^2
\end{eqnarray*}

\begin{eqnarray*}
(r^{k\lbar} -{r^k r^{\lbar}\over |\d r|_r^2})\TD r_{k\lbar} 
&=&4C(r^{k\kbar} -{r^k r^{\kbar}\over |\d r|_r^2})\TD |z_k|^2=4C(r^{k\kbar} -{r^k r^{\kbar}\over |\d r|_r^2})^2\\
%&=& 4C\delta_{k\lbar} ({|\d r|_r^2(1+4C |z_k|^2) -|r_k|^2 \over (1+4C |z_k|^2)^2|\d r|_r^2})
 \end{eqnarray*}
$$
\TD r_k=4C(r^{k\kbar}-{r^k r^{\kbar}\over |\d r|^2})\zbar_k
$$
and
$$
r_{\kbar}=(1+2C|z_k|^2) z_k +2\alpha \zbar_k.
$$
Thus by (5.3)
\begin{eqnarray*}
\Re r^k \TD r_k&=&4C\Re (r^{k\kbar}-{r^k r^{\kbar}\over |\d r|^2}) r^k \zbar_k\\
&\le &4C(r^{k\kbar}-{r^k r^{\kbar}\over |\d r|^2})r^{k\kbar}  (1+2\alpha+2C|z_k|^2)|z_k|^2\\
&=& {4C|z_k|^2 (1+2\alpha +2C|z_k|^2) \over (1+4C|z_k|^2)^2} \\
&\le&{2\alpha+1\over 8}
\end{eqnarray*}

\begin{eqnarray*}
(r^{k\lbar}-{r^k r^{\lbar}\over |\d r|^2}) \TD r_k \TD r_{\lbar}
&=&16 C^2 (r^{k\lbar}-{r^k r^{\lbar}\over |\d r|^2})\zbar_k z_{\ell} (r^{k\kbar}-{r^k r^{\kbar}\over |\d r|^2})(r^{\ell\lbar}-{r^\ell r^{\lbar}\over |\d r|^2})\\
&\le &16 C^2 r^{k\lbar}\zbar_k z_{\ell} (r^{k\kbar}-{r^k r^{\kbar}\over |\d r|^2})(r^{\ell\lbar}-{r^\ell r^{\lbar}\over |\d r|^2})\\
&\le& 4C {4C|z_k|^2 \over 1+4C|z_k|^2}(r^{k\kbar}-{r^k r^{\kbar}\over |\d r|^2})^2
\end{eqnarray*}
and
\begin{eqnarray*}
(r^{k\lbar}-{r^k r^{\lbar}\over |\d r|^2})(r^{i\qbar}-{r^i r^{\qbar}\over |\d r|^2})r^{p\jbar} r_{i\jbar k} r_{p\qbar\lbar}
&=&16C^2
(r^{k\lbar}-{r^k r^{\lbar}\over |\d r|^2})^2 r^{\ell \kbar} \zbar_k \delta_{ik} \delta_{jk} z_{\ell} \delta_{p\ell} \delta_{q\ell}\\
&=&16C^2 |z_k|^2
(r^{k\kbar}-{r^k r^{\kbar}\over |\d r|^2})^2 r^{k \kbar} \\
&=&4C  {4C |z_k|^2 \over (1+4C|z_k|^2)} (r^{k\kbar}-{r^{kk} |r_k|^2\over |\d r|^2})^2\\
\end{eqnarray*}
Therefore, since (5.1), we have
\begin{eqnarray*}
\lefteqn{(r^{k\lbar} -{r^k r^{\lbar}\over |\d r|_r^2})\TD r_{k\lbar} -(r^{k\lbar}-{r^k r^{\lbar}\over |\d r|^2})(r^{i\qbar}-{r^i r^{\qbar}\over |\d r|^2})r^{p\jbar} r_{i\jbar k} r_{p\qbar\lbar} 
-4C {4C|z_k|^2 \over 1+4C|z_k|^2}(r^{k\kbar}-{r^k r^{\kbar}\over |\d r|^2})^2
}\\
&=&4C(r^{k\kbar} -{r^k r^{\kbar}\over |\d r|_r^2})^2-4C  {4C |z_k|^2 \over (1+4C|z_k|^2)} (r^{k\kbar}-{r^{kk} |r_k|^2\over |\d r|^2})^2-
4C {4C|z_k|^2 \over 1+4C|z_k|^2}(r^{k\kbar}-{r^k r^{\kbar}\over |\d r|^2})^2
\\
&=&4C(1-2{4C|z_k|^2\over   1+4C|z_k|^2})(r^{k\kbar}-{r^k r^{\kbar} \over |\d r|^2})^2\\
&\ge& 0.
 \end{eqnarray*}
  Therefore,
  $$
 \tilde{E}(r)\ge- {2\Re r^k \TD r_k\over n+1}-a^{k\lbar}[r] {r^i r_{ik} r^{\jbar} r_{\jbar\lbar}\over (n+1)|\d r|^2}
 \ge -{(1+2\alpha)\over 4(n+1)}-{(2\alpha+1)^2\over 256(n+1)}
 $$
\begin{eqnarray*}
\lefteqn{1-{2\over n+1} \Re {r^i r^k r_{ik} \over |\d r|^2} +\tilde{E}(r) }\\
&\ge &1-{2\over n+1}\Re  {r^i r^i (\alpha+2C\zbar_i^2 )\over |\d r|^2} -{(1+2\alpha)\over 4(n+1)}-{(2\alpha+1)^2\over 256(n+1)} \\
&= &1-{2 \over n+1}\Re  { r^ {i\ibar} r_i^2 r^{i\ibar} (\alpha+2C\zbar_i^2 )\over |\d r|^2}  -{(1+2\alpha)\over 4(n+1)}-{(2\alpha+1)^2\over 256 (n+1)}\\
&\ge  &1-{2\alpha  \over n+1}   -{(1+2\alpha)\over 4(n+1)}-{(2\alpha+1)^2 \over 256(n+1)}\\
&>&1-{10\alpha+1\over 4(n+1)}-{10 \over 256(n+1)}\\
&\ge & 1-{23\over 24}-{1\over 25}\\
&>&0
\end{eqnarray*}
if $n\ge 2$ and $\alpha\le 21/20$.
  Therefore,  $ D$ is strictly super-pseudoconvex and the proof is complete.\epf

\bigskip

\noindent Mailing address: 

Department of Mathematics, University of California, 
  Irvine, CA 92697--3875.

\noindent E-mail: \quad sli@math.uci.edu
\end{document}